\documentclass[a4paper,12pt]{article}
\usepackage{comment}
\usepackage{cite}
\usepackage{amsmath,amssymb,amsfonts,amsthm,mathrsfs,empheq}
\usepackage[T1]{fontenc}
\usepackage[utf8]{inputenc}
\usepackage{graphicx}
\usepackage{fancyhdr}
\usepackage{float}
\usepackage{xcolor}
\usepackage{authblk}
\usepackage[hyphens]{url}
\usepackage{hyperref}
\usepackage[]{breakurl}
\usepackage{tikz}
\usetikzlibrary{decorations.markings,arrows.meta,positioning}
\usepackage{geometry}
\theoremstyle{plain}

\newtheorem{definition}{Definition}
\newtheorem{remark}{Remark}
\newtheorem{example}{Example}

\title{Eulerian polynomials and the alternating sum of excedances}

\author[$\dagger$]{Jean-Christophe {\sc Pain}$^{1,2,}$\footnote{jean-christophe.pain@cea.fr}\\
\small
$^1$CEA, DAM, DIF, F-91297 Arpajon, France\\
$^2$Universit\'e Paris-Saclay, CEA, Laboratory Matière en Conditions Extr\^emes,\\ 
F-91680 Bruy\`eres-le-Ch\^atel, France
}

\date{}

\begin{document}

\maketitle

\begin{abstract}
Tangent numbers $T_{2n-1}$, which enumerate alternating permutations of odd length, play a prominent role in the Taylor series expansion of the tangent function $\tan(x)$. In this work, we adopt a combinatorial approach based on the excedance statistic of permutations, which allows us to interpret the coefficients of the tangent series in a structural and enumerative way. Using this framework, we establish a classical identity that relates the alternating sum of excedances to the hyperbolic tangent function. This perspective highlights deep connections with Eulerian polynomials, provides a combinatorial interpretation of tangent numbers, and links these sequences to Genocchi numbers and related arithmetic properties. The approach not only unifies analytic and combinatorial viewpoints but also opens the way to generalizations to other permutation statistics and families of specialized permutations.
\end{abstract}

\section{Introduction}

Combinatorics and the study of special sequences in mathematics often reveal deep and surprising connections between seemingly unrelated structures. Among these, the tangent numbers $T_{2n-1}$, which enumerate alternating permutations, play a major role. Initially appearing in the classical Taylor series expansion of the tangent function,
\[
\tan x = \sum_{n=1}^{\infty} \frac{T_{2n-1}}{(2n-1)!} x^{2n-1},
\]
these numbers have a rich combinatorial interpretation: they enumerate the alternating permutations of odd length, that is, permutations of $\{1,2,\dots,2n-1\}$ in which entries rise and fall in an alternating pattern (up-down-up-...).

\begin{figure}[!ht]
\centering
\begin{tikzpicture}[scale=1, every node/.style={circle, draw, minimum size=8mm, inner sep=1mm}, 
    >=Stealth, node distance=1.5cm]

% Nodes
\node (1) {1};
\node (2) [right=of 1] {2};
\node (3) [right=of 2] {3};
\node (4) [right=of 3] {4};
\node (5) [right=of 4] {5};

% Example alternating permutation: 2 4 1 5 3
\draw[->, thick, blue] (1) to[out=45,in=135] (2);
\draw[->, thick, red] (2) to[out=-45,in=135] (4);
\draw[->, thick, blue] (3) to[out=45,in=-135] (1);
\draw[->, thick, red] (4) to[out=45,in=-135] (5);
\draw[->, thick, blue] (5) to[out=-45,in=-135] (3);

% Legend
\node[draw=none, fill=none, right=1cm of 5] (leg1) {\textcolor{blue}{\textbf{Up}}};
\node[draw=none, fill=none, below=0.4cm of leg1] (leg2) {\textcolor{red}{\textbf{Down}}};

\end{tikzpicture}
\caption{Example of an alternating permutation for $n=5$, corresponding to a term counted by $T_5$. Blue arrows indicate ``up'' steps and red arrows indicate ``down'' steps.}\label{fig:alternating_permutation}
\end{figure}
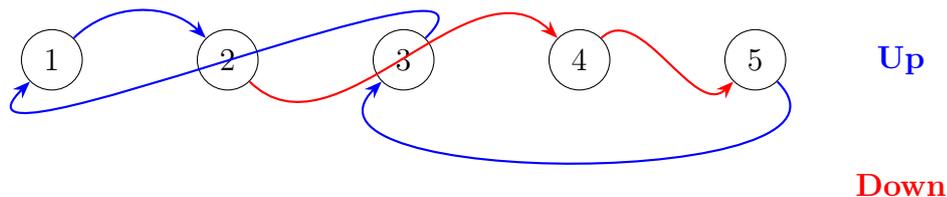

The study of alternating permutations dates back to early combinatorial investigations by Andr\'e in the 19th century, and they have since been connected to a variety of mathematical areas, including the analysis of continued fractions, Bernoulli and Genocchi numbers, and the arithmetic properties of special sequences. Beyond their historical significance, tangent numbers and alternating permutations appear naturally in enumerative combinatorics, algebraic combinatorics, and even in probabilistic models, highlighting the ubiquity of these structures in mathematics. Figure~\ref{fig:alternating_permutation} illustrates an alternating permutation of length $n=5$, with blue and red arrows representing ``up'' and ``down'' steps, highlighting the zigzag structure counted by tangent numbers.

% TikZ figure

In parallel, permutation statistics, such as the excedance number, offer a powerful framework for understanding the structural properties of permutations. For a permutation $\sigma \in \mathfrak{S}_n$, the number of excedances is defined by $\#\{ i \in \{1,\dots,n\} \;|\; \sigma(i) > i \}$. This statistic captures the ``upward steps'' in a permutation and generalizes classical notions of descents and ascents. Studying the distribution of excedances over all permutations of a given length provides a bridge between combinatorial enumeration and analytic techniques, such as generating functions. It is worth mentioning that Various $q$‑enumerations of alternating permutations, including those studied by Josuat‑Vergès \cite{JosuatVerges}, illustrate the rich combinatorial structures associated with the excedance statistic.

In this work, we adopt this combinatorial perspective based on excedances to revisit the tangent numbers. Using generating functions of Eulerian polynomials and the alternating sum of excedances, we recover a classical identity from enumerative combinatorics: evaluating the exponential generating function of Eulerian polynomials at $t=-1$ yields the hyperbolic tangent function, which encodes the tangent numbers. This approach provides a standard combinatorial interpretation of the alternating sum of excedances in terms of $T_{2n-1}$. This work not only recovers classical results about tangent numbers and alternating permutations but also interprets connections with Genocchi numbers and provides a conceptual explanation for parity and integrality properties.

In Section 2, we recall fundamental definitions concerning permutations, the excedance statistic, and Eulerian polynomials, establishing the notations and generating functions that will be used throughout the paper. Section 3 presents an identity for the alternating sum of excedances and demonstrates its direct connection with the hyperbolic tangent function and tangent numbers. In Section 4, we explore several derived properties, including parity results, arithmetic consequences, and links with Genocchi numbers.

\section{Fundamental definitions}

\subsection{Permutations and excedances}

A permutation of $n$ elements is a bijection $\sigma: \{1,2,\dots,n\} \to \{1,2,\dots,n\}$. Permutations can be analyzed through various statistics that capture their combinatorial structure. One of the most classical and widely studied statistics is the excedance, which measures the number of positions $i$ in a permutation where the value exceeds the position index.

\begin{definition}[Excedance]
For $\sigma \in \mathfrak{S}_n$, the number of excedances $exc(\sigma)$ is defined as
\[
exc(\sigma) := \#\{ i \in \{1,\dots,n\} : \sigma(i) > i \}.
\]
\end{definition}

\begin{remark}
The excedance statistic is closely related to other classical permutation statistics, such as descents and ascents, and plays a central role in enumerative combinatorics. It encodes information about the ``upward steps'' in a permutation, and appears naturally in the study of Eulerian polynomials, alternating permutations, and tangent numbers.
\end{remark}

\begin{example}
Consider the permutation $\sigma = (2,4,1,3) \in \mathfrak{S}_4$. The excedances are calculated by comparing each value to its position:
\[
\sigma(1)=2>1, \quad \sigma(2)=4>2, \quad \sigma(3)=1\not>3, \quad \sigma(4)=3\not>4.
\]
Thus, $exc(\sigma) = 2$. This means that two elements of the permutation exceed their positions.
\end{example}

\begin{remark}
Understanding the distribution of excedances across all permutations of length $n$ provides combinatorial insight into the structure of $\mathfrak{S}_n$, and sets the stage for analytic tools, such as generating functions, that will be used to derive identities connecting permutation statistics to classical sequences like tangent numbers.
\end{remark}

\subsection{Eulerian polynomials}

The Eulerian polynomials $A_n(t)$ provide a classical and powerful way to encode the excedance statistic over the symmetric group $\mathfrak{S}_n$:
\[
A_n(t) = \sum_{\sigma \in \mathfrak{S}_n} t^{exc(\sigma)+1}, \quad \text{with the $+1$ included by convention so that } A_n(1) = n!.
\]
Eulerian polynomials have a rich combinatorial structure and satisfy numerous interesting properties \cite{Stanley}. For instance, the coefficient of $t^k$ in $A_n(t)$ counts the number of permutations in $\mathfrak{S}_n$ with exactly $k-1$ excedances, and they are closely related to other classical permutation statistics such as descents, ascents, and runs. Historically, they were first studied by D. Andr\'e \cite{Andre1879,Andre1881} and later extensively analyzed in enumerative combinatorics, serving as a central example of how generating functions encode combinatorial information (see also Foata and Sch\"utzenberger \cite{FoataSchutzenberger}).

The exponential generating function (EGF) of Eulerian polynomials is given by the classical identity:
\begin{equation}
\Phi(x,t) := \sum_{n=0}^{\infty} A_n(t) \frac{x^n}{n!} = \frac{t-1}{t - e^{x(t-1)}}.
\label{eq:FGE}
\end{equation}

\begin{remark}
This EGF not only provides an analytic tool to study the distribution of excedances but also allows one to connect permutation statistics with classical functions. In particular, evaluating $\Phi(x,t)$ at specific values of $t$ leads to hyperbolic functions and sequences such as the tangent numbers, thereby illustrating a deep interplay between combinatorics and analysis. Furthermore, the EGF encodes recurrences, parity properties, and identities for Eulerian numbers in a compact and elegant form.
\end{remark}

\section{An identity for the alternating sum of excedances}

In this section, we establish a classical identity connecting the alternating sum of excedances in permutations to tangent numbers and the hyperbolic tangent function. We provide a coherent derivation using generating functions and highlight the combinatorial interpretation.  

\subsection{Definition of the alternating sum}

Let us define the alternating sum of excedances for permutations of length \(n\) as  
\[
S_n := \sum_{\sigma \in \mathfrak{S}_n} (-1)^{exc(\sigma)}.
\]
This sum encodes the difference between the number of permutations with an even number of excedances and those with an odd number. Our goal is to express \(S_n\) explicitly and to connect it with classical sequences, notably the tangent numbers \(T_{2n-1}\).  

\subsection{Exponential generating function approach}

From the definition of Eulerian polynomials \(A_n(t) = \sum_{\sigma \in \mathfrak{S}_n} t^{exc(\sigma)+1}\), we immediately see that  
\[
A_n(-1) = \sum_{\sigma \in \mathfrak{S}_n} (-1)^{exc(\sigma)+1} = - S_n.
\]
Thus, up to a global sign, the alternating sum \(S_n\) is encoded in the exponential generating function (EGF) of Eulerian polynomials evaluated at \(t=-1\):
\[
\Phi(x,-1) := \sum_{n=0}^{\infty} A_n(-1) \frac{x^n}{n!} = - \sum_{n=0}^{\infty} S_n \frac{x^n}{n!}.
\]
Substituting \(t=-1\) into the classical EGF formula  
\[
\Phi(x,t) = \frac{t-1}{t - e^{x(t-1)}},
\]
we obtain
\[
\Phi(x,-1) = \frac{-2}{-1 - e^{-2x}} = \frac{2}{1 + e^{-2x}}.
\]
Multiplying numerator and denominator by \(e^x\), this simplifies to
\[
\Phi(x,-1) = \frac{2 e^x}{e^x + e^{-x}} = 1 + \tanh x.
\]
Equating the two expressions of \(\Phi(x,-1)\) yields
\[
- \sum_{n=0}^{\infty} S_n \frac{x^n}{n!} = 1 + \sum_{n=1}^{\infty} \frac{T_{2n-1}}{(2n-1)!} x^{2n-1},
\]
where \(T_{2n-1}\) denotes the \((2n-1)\)-th tangent number.  

\subsection{Explicit sum rule}

From the coefficient comparison, we immediately obtain  
\[
S_0 = 1, \quad S_{2n} = 0, \quad S_{2n-1} = (-1)^{\,n-1} T_{2n-1}, \quad n \ge 1.
\]
This explicit formula has a natural combinatorial interpretation. For even \(n\), each permutation with an even number of excedances is exactly canceled by a permutation with an odd number, resulting in \(S_{2n}=0\). This cancellation corresponds to the fact that \(\tanh x\) is an odd function, so only odd powers appear in its Taylor expansion. For odd \(n = 2m-1\), the alternating sum does not vanish; the contributions combine to produce a signed count \(S_{2m-1} = (-1)^{m-1} T_{2m-1}\), directly linking the alternating sum of excedances to the tangent numbers.

\subsection{Connection with Bernoulli numbers}

Tangent numbers can also be expressed explicitly in terms of Bernoulli numbers $B_{2n}$:
\begin{equation}
T_{2n-1} = (-1)^{\,n-1} \frac{2^{2n}\,(2^{2n}-1)}{2n} \, B_{2n}.
\end{equation}
This explicit formula confirms that the sign $(-1)^{exc(\sigma)}$ in the alternating sum matches the classical sign convention for tangent numbers derived from Bernoulli numbers. Combining the EGF approach with this expression provides a clear analytic-combinatorial bridge: the coefficients of $x^{2n-1}/(2n-1)!$ in $\tanh x$ correspond exactly to $T_{2n-1}/(2n-1)!$, which enumerate alternating permutations of odd length, weighted by $(-1)^{n-1}$.

\subsection{Combinatorial interpretation}

This perspective provides an intuitive understanding of why the odd part of the generating function corresponds to tangent numbers, and why the sum vanishes for even lengths. We will see in Section~4 how these combinatorial interpretations lead to arithmetic and recurrence properties.

\section{Derived properties}

The identity relating alternating excedance sums to tangent numbers leads naturally to several interesting combinatorial and arithmetic properties.

\subsection{Integrality} 

The tangent numbers $T_{2n-1}$ are integers, and this integrality is reflected combinatorially in the fact that the alternating sum of excedances over $\mathfrak{S}_{2n-1}$ counts signed combinatorial objects (see also Ref. \cite{Pain2026} for further analysis of their arithmetic structure). More precisely, each permutation contributes $\pm 1$, so the total sum is necessarily an integer. This also explains why the coefficients of the series expansion of $\tanh x$ can be expressed in terms of $T_{2n-1}/(2n-1)!$ even though these coefficients are rational when written in standard factorial form.

\subsection{Parity} 

For even $n>0$, the alternating sum vanishes:
\[
\sum_{\sigma \in \mathfrak{S}_n} (-1)^{exc(\sigma)} = 0.
\]
This follows directly from the fact that $\tanh x$ is an odd function: its Taylor expansion contains only odd powers of $x$. Combinatorially, this reflects a perfect cancellation between permutations with an even and odd number of excedances.

\subsection{Connection with Genocchi numbers} 

The alternating sums are closely related to the Genocchi numbers $G_{n+1}$:
\[
S_n = \sum_{\sigma \in \mathfrak{S}_n} (-1)^{exc(\sigma)} = (-1)^{\lfloor (n+1)/2 \rfloor} G_{n+1}.
\]
Genocchi numbers, which can be defined via the generating function
\[
\frac{2t}{e^t+1} = \sum_{n=1}^{\infty} G_n \frac{t^n}{n!},
\]
appear naturally in the study of alternating permutations and provide a bridge between tangent numbers and Bernoulli numbers. This connection also reveals a rich arithmetic structure underlying the excedance statistic, as the Genocchi numbers encode divisibility properties and congruences reminiscent of those of tangent numbers. Furthermore, we can exploit known arithmetic properties of $G_n$ to deduce corresponding results for $S_n$. For instance, $S_{2n} = 0$, reflecting the vanishing of $G_{2n+1}$ for $n\ge 1$, and
$S_{2n-1}$ is divisible by $2$, since $G_{2n}$ is divisible by $2$. We have also some congruences modulo $4$: $S_{4n-1} \equiv 0 \pmod{4}$, $S_{4n+1} \equiv 2 \pmod{4}$, etc. The Genocchi numbers satisfy classical recurrences of the form
\[
G_n = - \sum_{k=1}^{n-1} \binom{n}{k} G_k, \quad n\ge 2, \quad G_1=1.
\]
Through the correspondence $S_n = (-1)^{\lfloor (n+1)/2 \rfloor} G_{n+1}$, this immediately induces recurrences for $S_n$ up to a predictable sign factor:
\[
S_n = \sum_{k=1}^{n-1} (-1)^{f(n,k)} \binom{n+1}{k} S_k,
\]
where the exponent $f(n,k)$ keeps track of the alternating signs arising from the mapping between $G_{n+1}$ and $S_n$:
\[
f(n,k) := \left\lfloor \frac{n+1}{2} \right\rfloor - \left\lfloor \frac{k+1}{2} \right\rfloor.
\]
Combinatorially, these recurrences can be interpreted as \emph{insertion rules}: one can construct all permutations of length $n$ contributing to $S_n$ by inserting the largest element into smaller permutations of length $k$ in a way that preserves the alternating excedance structure. This provides an elegant explanation of why the sums vanish for even $n$, why the signs alternate for odd $n$, and how the arithmetic structure of $S_n$ is inherited from the classical identities of Genocchi numbers.

\subsection{Recurrence relations}

The alternating sums $S_n = \sum_{\sigma \in \mathfrak{S}_n} (-1)^{exc(\sigma)}$ satisfy simple parity-based relations:
\[
S_0 = 1, \quad S_{2n} = 0, \quad S_{2n-1} = (-1)^{n-1} T_{2n-1}, \quad n \ge 1.
\]
One can also express $S_{n+1}$ recursively in terms of smaller sums using combinatorial insertion of the largest element:
\[
S_{n+1} = \sum_{k=0}^{n} (-1)^k \binom{n}{k} S_k,
\]
but one must take care with the sign convention to ensure consistency with $S_{2n} = 0$. This reflects how inserting the largest element into permutations of length $n$ affects the number of excedances. This viewpoint provides a unified combinatorial explanation for the cancellations and parity phenomena observed in the alternating sums, linking signed permutations directly to tangent and Genocchi numbers.

\section{Conclusion}

The excedance statistic provides a powerful and versatile framework for understanding the arithmetic and combinatorial structure of tangent numbers and Eulerian polynomials. The identity relating the alternating sum of excedances to the hyperbolic tangent function not only offers a direct analytic link to classical sequences such as tangent numbers but also highlights deep structural properties of permutations, including parity, signed counts, and cancellations. 

Moreover, this approach naturally extends to related combinatorial families, such as Genocchi numbers and alternating permutations of even length, revealing intricate connections between generating functions, permutation statistics, and classical number sequences. Beyond these concrete examples, the methods presented here suggest a broader perspective: similar analytic-combinatorial techniques may be applied to study other permutation statistics, weighted sums, or generalized Eulerian polynomials, potentially uncovering new identities and congruences in algebraic combinatorics. 

In summary, the interplay between generating functions, excedance statistics, and classical sequences underscores the rich synergy between combinatorial theory and analytic methods, opening avenues for further exploration in both enumerative combinatorics and the arithmetic of special sequences.

\end{document}